# MATCHING WITH SHIFT FOR ONE-DIMENSIONAL GIBBS MEASURES

By P. Collet, C. Giardina and F. Redig

*CNRS UMR 7644, Eindhoven University and Universiteit Leiden*

We consider matching with shifts for Gibbsian sequences. We prove that the maximal overlap behaves as $c \log n$, where $c$ is explicitly identified in terms of the thermodynamic quantities (pressure) of the underlying potential. Our approach is based on the analysis of the first and second moment of the number of overlaps of a given size. We treat both the case of equal sequences (and nonzero shifts) and independent sequences.

**1. Introduction.** In sequence alignment one wants to detect significant similarities between two (e.g., genetic or protein) sequences. In order to distinguish "significant" similarities, one has to compute the probability that a similarity of a certain size occurs for two independent sequences. The symbols in the sequences are, however, not necessarily occurring independently. From the point of view of statistical mechanics, it is quite natural to assume that the symbols in the sequence are generated according to a stationary Gibbs measure: this is the equilibrium measure which maximizes the entropy under physical constraints such as energy conservation. A priori there is no reason to assume that the symbols (bases) in, for example, a DNA sequence, are i.i.d. or even Markov. It can, however, be plausible to assume that there is an underlying Markov chain of which the symbol sequence is a reduction: in that case we arrive at a so-called hidden Markov chain, and it is well known that hidden Markov chains have generically infinite memory (though the symbol at a particular location only exponentially weakly depends on symbols far away). Therefore, proposing a Gibbs measure with exponentially decaying interaction as a model for the sequence seems quite natural. Besides motivation coming from sequence alignment, also in dynamical systems, [4] one can ask for the probability of having a large "overlap" in a trajectory of length $n$, but without specifying the location of the piece of trajectory that is repeated. It is clear that this probability is related to the entropy, but not in such a straightforward way as the return time. In (hyperbolic) dynamical systems, by coding and partitioning, one again naturally arrives at Gibbs measures with exponentially decaying interactions.

The first nontrivial problem associated with sequence alignment is the comparison of two sequences where it is allowed to shift one sequence w.r.t. the other. Remark that this problem is not easy even in the case of independent symbols in the sequence, because one allows for shifting one sequence w.r.t. the other. The comparison consists in the simplest case in finding the maximal number of consecutive equal symbols. Given two (independent) i.i.d. sequences, in [5] and [6] it is proved that the maximal overlap, allowing shifts, behaves for large sequence length as $c \log n + X$, where $n$ is the length of both sequences, $c$ is a constant depending on the distribution of the sequence, and where $X$ is a random variable with a Gumbel distribution. The fact that $c \log(n)$ is the good scale can be easily understood intuitively: it corresponds to the maximum of order $n$ weakly dependent variables. However, even in the case of i.i.d. sequences, it is not so easy to make that intuition rigorous, as we allow shifts. In fact, the results of [5] and [6] are based on large deviations, together with an analysis of random walk excursions. As the proofs use a form of permutation invariance, they cannot be extended to non-i.i.d. cases. In [9] the maximal alignment with shift is shown for Markov sequences, which requires a theory of excursions of random walk with Markovian increments.

In this paper we focus on the more elementary question of showing that the maximal overlap allowing shifts behaves as $c \log n$, but now in the context of general Gibbsian sequences. We also allow to match a sequence with *itself* (where of course we have to restrict to nonzero shifts). The constant $c$ is explicitly identified and related to thermodynamic quantities associated to the potential of the underlying Gibbs measure.

Our approach is based on a first and second moment analysis of the random variable $N(\sigma, n, k)$ that counts the number of shift-matches of size $k$ in a sequence $\sigma$ of length $n$. One easily identifies the scale $k = k_n = c \log(n)$ which discriminates the region where the first moment $\mathbb{E} N(\sigma, n, k_n)$ goes to zero (as $n \to \infty$) from the region where $\mathbb{E} N(\sigma, n, k)$ diverges. Via a second moment estimate, we then prove that this scale also separates the $N(\sigma, n, k) \to 0$ versus $N(\sigma, n, k) \to \infty$ (convergence in probability) region.

Our paper is organized as follows: in Section 2 we introduce the basic preliminaries about Gibbs measures, in Section 3 we analyze the first moment of $N$ in the case of matching a sequence with itself and in Section 4 we study the second moment. In Section 5 we treat the case of two independent (Gibbsian) sequences with the same and with different marginal distributions.

**2. Definitions and preliminaries.** We consider random stationary sequences [8] $\sigma = \{\sigma(i) : i \in \mathbb{Z}\}$ on the lattice $\mathbb{Z}$, where $\sigma(i)$ takes values in a finite set $\mathcal{A}$. The joint distribution of $\{\sigma(i) : i \in \mathbb{Z}\}$ is denoted by $\mathbb{P}$. We treat the case where $\mathbb{P}$ is a Gibbs measure with exponentially decaying interaction; see Section 2.3 below for details. The configuration space $\Omega = \mathcal{A}^\mathbb{Z}$ is endowed with the product topology (making it into a compact metric space). The set of finite subsets of $\mathbb{Z}$ is denoted by $\mathcal{S}$. For $V, W \in \mathcal{S}$, we put $d(V, W) = \min\{|i - j| : i \in V, j \in W\}$. For $V \in \mathcal{S}$, the diameter is defined via $\mathrm{diam}(V) = \max\{|i - j|, i, j \in V\}$. For $V \in \mathcal{S}$, $\mathcal{F}_V$ is the sigma-field generated by $\{\sigma(i) : i \in A\}$. For $V \in \mathcal{S}$, we put $\Omega_V = \mathcal{A}^V$. For $\sigma \in \Omega$ and $V \in \mathcal{S}$, $\sigma_V \in \Omega_V$ denotes the restriction of $\sigma$ to $V$. For $i \in \mathbb{Z}$ and $\sigma \in \Omega$, $\tau_i \sigma$ denotes the translation of $\sigma$ by $i : \tau_i \sigma(j) = \sigma(i + j)$. For a local event $E \subseteq \Omega$, the dependence set of $E$ is defined by the minimal $V \in \mathcal{S}$ such that $E$ is $\mathcal{F}_V$ measurable. We denote $\mathbb{1}$ for the indicator function.

2.1. *Patterns and cylinders.* For $n \in \mathbb{N}, n \geq 1$, let $C_n = [1, n] \cap \mathbb{Z}$. An element $A_n \in \Omega_{C_n}$ is called a $n$-pattern or a pattern of size $n$. For a pattern $A_n \in \Omega_{C_n}$, we define the corresponding cylinder $\mathscr{C}(A_n) = \{\sigma \in \Omega : \sigma_{C_n} = A_n\}$. The collection of all $n$-cylinders is denoted by $\mathscr{C}_n = \bigcup_{A_n \in \Omega_{C_n}} \mathscr{C}(A_n)$. Sometimes, to denote the probability of the cylinder associated to the pattern $A_n$, we will use the abbreviation

(2.1) $$\mathbb{P}(A_n) := \mathbb{P}(\mathscr{C}(A_n)) = \mathbb{P}(\sigma_{C_n} = A_n).$$

For $A_k = (\sigma(1), \sigma(2), \ldots, \sigma(k))$ a $k$-pattern and $1 \leq i \leq j \leq n$, we define the pattern $A_k(i, j)$ to be the pattern of length $j - i + 1$ consisting of the symbols $(\sigma(i), \sigma(i+1), \ldots, \sigma(j))$. For two patterns $A_k, B_l$, we define their concatenation $A_k B_l$ to be the pattern of length $k + l$ consisting of the $k$ symbols of $A_k$ followed by the $l$ symbols of $B_l$. Concatenation of three or more patterns follows obviously from this.

2.2. *Shift-matches.* We will study properties of the following basic quantities.

DEFINITION 2.1 (Number of shift-matches). For every configuration $\sigma \in \Omega$ and for every $n \in \mathbb{N}$, $k \in \mathbb{N}$, with $k \leq n$, we define the number of matches with shift of length $k$ up to $n$ as

$$N(\sigma, n, k) = \frac{1}{2} \sum_{i=0}^{n-k} \sum_{j=0, j\neq i}^{n-k} \mathbb{1}\{(\tau_i \sigma)_{C_k} = (\tau_j \sigma)_{C_k}\}$$

(2.2) $$= \sum_{i \neq j = 0}^{n-k} \mathbb{1}(\sigma(i+1) = \sigma(j+1), \sigma(i+2) = \sigma(j+2), \ldots,$$

$$\sigma(i+k) = \sigma(j+k)).$$

DEFINITION 2.2 (Maximal shift-matching). For every configuration $\sigma \in \Omega$ and for every $n \in \mathbb{N}$, we define $M(\sigma, n)$ to be the maximal length of a shift-matching up to $n$, that is the maximal $k \in \mathbb{N}$ (with $k \leq n$) such that there exist $i \in \mathbb{N}$ and $j \in \mathbb{N}$ (with $0 \leq i < j \leq n - k$) satisfying

$$(2.3) \qquad (\tau_i \sigma)_{C_k} = (\tau_j \sigma)_{C_k},$$

where we adopt the convention $\max(\varnothing) = 0$.

DEFINITION 2.3 (First occurrence of a shift-matching). For every configuration $\sigma \in \Omega$ and for every $k \in \mathbb{N}$, we define $T(\sigma, k)$ to be the first occurrence of a shift-match, that is, the minimal $n \in \mathbb{N}$ (with $k \leq n$) such that there exist $i \in \mathbb{N}$ and $j \in \mathbb{N}$ (with $0 \leq i < j \leq n - k$) satisfying

$$(2.4) \qquad (\tau_i \sigma)_{C_k} = (\tau_j \sigma)_{C_k},$$

where we adopt the convention $\min(\varnothing) = \infty$.

The following proposition follows immediately from these definitions.

PROPOSITION 2.4. *The probability distributions of the previous quantities are related by the following "duality" relations:*

$$(2.5) \qquad \mathbb{P}(N(\sigma, n, k) = 0) = \mathbb{P}(M(\sigma, n) < k) = \mathbb{P}(T(\sigma, k) > n).$$

2.3. *Gibbs measures.* We now state our assumptions on $\mathbb{P}$, and recall some basic facts about Gibbs measures [11]. The reader familiar with this can skip this section.

We choose for $\mathbb{P}$ the unique Gibbs measure corresponding to an exponentially decaying translation-invariant interaction. In dynamical systems language this corresponds to the unique equilibrium measure of a Hölder continuous potential.

2.3.1. *Interactions.*

DEFINITION 2.5. A translation-invariant interaction is a map

$$(2.6) \qquad U : \mathcal{S} \times \Omega \to \mathbb{R},$$

such that the following conditions are satisfied:

1. For all $A \in \mathcal{S}$, $\sigma \mapsto U(A, \sigma)$ is $\mathcal{F}_A$-measurable.
2. *Translation invariance*:

$$(2.7) \qquad U(A + i, \tau_{-i} \sigma) = U(A, \sigma) \qquad \forall A \in \mathcal{S},\ i \in \mathbb{Z},\ \sigma \in \Omega.$$

3. *Exponential decay*: there exist $\gamma > 0$ such that

$$\|U\|_\gamma := \sum_{A \ni 0} e^{\gamma \operatorname{diam}(A)} \sup_{\sigma \in \Omega} |U(A,\sigma)| < \infty. \tag{2.8}$$

The set of all such interactions is denoted by $\mathcal{U}$. Here are some standard examples of elements of $\mathcal{U}$:

1. Ising model with magnetic field $h$: $\mathcal{A} = \{-1, 1\}$, $U(\{i, i+1\}, \sigma) = J\sigma_i\sigma_{i+1}$, $U(\{i\}, \sigma) = h\sigma_i$ and all other $U(A, \sigma) = 0$. Here $J, h \in \mathbb{R}$. If $J < 0$, we have the standard ferromagnetic Ising model.
2. General finite range interactions. An interaction $U$ is called *finite-range* if there exists an $R > 0$ such that $U(A, \sigma) = 0$ for all $A \in \mathcal{S}$ with $\operatorname{diam}(A) > R$.
3. Long range Ising models $U(\{i, j\}, \sigma) = J_{j-i}\sigma_i\sigma_j$ with $|J_k| \leq e^{-\gamma k}$ for some $\gamma > 0$ and $U(A, \sigma) = 0$ for all other $A \in \mathcal{S}$.

2.3.2. *Hamiltonians.* For $U \in \mathcal{U}$, $\zeta \in \Omega$, $\Lambda \in \mathcal{S}$, we define the finite-volume Hamiltonian with boundary condition $\zeta$ as

$$H_\Lambda^\zeta(\sigma) = \sum_{A \cap \Lambda \neq \varnothing} U(A, \sigma_\Lambda \zeta_{\Lambda^c}) \tag{2.9}$$

and the Hamiltonian with free boundary condition as

$$H_\Lambda(\sigma) = \sum_{A \subseteq \Lambda} U(A, \sigma), \tag{2.10}$$

which depends only on the spins inside $\Lambda$. In particular, for $A_k$ a pattern, $\sigma \in \mathscr{C}(A_k)$, $H_{C_k}(\sigma)$ depends only on $A_k$. We will denote, therefore,

$$H(\mathscr{C}(A_k)) = H_{C_k}(\sigma)$$

for $\sigma \in \mathscr{C}(A_k)$.

Corresponding to the Hamiltonian in (2.9), we have the finite-volume Gibbs measures $\mathbb{P}_\Lambda^{U,\zeta}$, $\Lambda \in \mathcal{S}$, defined on $\Omega$ by

$$\int f(\xi) \, d\mathbb{P}_\Lambda^{U,\zeta}(\xi) = \sum_{\sigma_\Lambda \in \Omega_\Lambda} f(\sigma_\Lambda \zeta_{\Lambda^c}) \frac{e^{-H_\Lambda^\zeta(\sigma)}}{Z_\Lambda^\zeta}, \tag{2.11}$$

where $f$ is any continuous function and $Z_\Lambda^\zeta$ denotes the partition function normalizing $\mathbb{P}_\Lambda^{U,\zeta}$ to a probability measure:

$$Z_\Lambda^\zeta = \sum_{\sigma_\Lambda \in \Omega_\Lambda} e^{-H_\Lambda^\zeta(\sigma)}. \tag{2.12}$$

2.3.3. *Gibbs measures with given interaction.* For a probability measure $\mathbb{P}$ on $\Omega$, we denote by $\mathbb{P}_\Lambda^\zeta$ the conditional probability distribution of $\sigma(i), i \in \Lambda$, given $\sigma_{\Lambda^c} = \zeta_{\Lambda^c}$. Of course, this object is only defined on a set of $\mathbb{P}$-measure one. For $\Lambda \in \mathcal{S}, \Gamma \in \mathcal{S}$ and $\Lambda \subseteq \Gamma$, we denote by $\mathbb{P}_\Gamma(\sigma_\Lambda|\zeta)$ the conditional probability to find $\sigma_\Lambda$ inside $\Lambda$, given that $\zeta$ occurs in $\Gamma \setminus \Lambda$.

DEFINITION 2.6. For $U \in \mathcal{U}$, we call $\mathbb{P}$ a Gibbs measure with interaction $U$ if its conditional probabilities coincide with the ones prescribed in (2.11), that is, if

(2.13) $$\mathbb{P}_\Lambda^\zeta = \mathbb{P}_\Lambda^{U,\zeta} \qquad \mathbb{P}\text{-a.s. } \Lambda \in \mathcal{S}, \ \zeta \in \Omega.$$

In our situation, with $U \in \mathcal{U}$, the Gibbs measure $\mathbb{P}$ corresponding to $U$ is unique. Moreover, it satisfies the following strong mixing condition: for all $V, W \in \mathcal{S}$ and all events $A \in \mathcal{F}_V$, $B \in \mathcal{F}_W$,

(2.14) $$\left|\frac{\mathbb{P}(A \cap B)}{\mathbb{P}(B)} - \mathbb{P}(A)\right| \leq e^{-c\, d(V,W)},$$

where $c > 0$ depends of course on the interaction $U$.

2.4. *Thermodynamic quantities.* We now recall some definitions of basic important statistical mechanics quantities.

DEFINITION 2.7. The pressure $p(U)$ of the Gibbs measure $\mathbb{P}$ associated with the interaction $U$ is defined as

(2.15) $$p(U) = \lim_{n \to \infty} \frac{1}{n} \log Z_n,$$

where

$$Z_n = \sum_{\sigma_{C_n} \in \Omega_{C_n}} \exp\left(-\sum_{A \subseteq C_n} U(A, \sigma)\right)$$

is the partition function with the free boundary conditions.

DEFINITION 2.8. The entropy $s(U)$ of the Gibbs measure $\mathbb{P}$ associated with the interaction $U$ is defined as

(2.16) $$s(U) = \lim_{n \to \infty} -\frac{1}{n} \sum_{A_n \in \Omega_{C_n}} \mathbb{P}(\mathscr{C}(A_n)) \log \mathbb{P}(\mathscr{C}(A_n)).$$

In terms of the interaction $U$, we have the following basic thermodynamic relation between pressure, entropy and the Gibbs measure $\mathbb{P}$ corresponding to $U$:

(2.17) $$s(U) = p(U) + \int f_U \, d\mathbb{P},$$

where
$$f_U(\sigma) = \sum_{A \ni 0} \frac{U(A, \sigma)}{|A|}$$

denotes the average internal energy per site.

We also have the following relation between $f_U$ and the Hamiltonian:

$$(2.18) \qquad H_\Lambda^\xi(\sigma) = \sum_{i \in \Lambda} \tau_i f_U(\sigma) + O(1),$$

where $O(1)$ is a quantity which is uniformly bounded in $\Lambda, \sigma, \xi$.

The function $f_U$ is what is called the potential in the dynamical systems literature. An exponentially decaying interaction $U$ then corresponds to a Hölder continuous potential $f_U$.

The following is a standard property of (one-dimensional) Gibbs measures with interaction $U \in \mathcal{U}$. For the proof, see [3], page 7. See also [7], pages 164–165 for properties of one-dimensional Gibbs measures.

PROPOSITION 2.9. *For the unique Gibbs measure $\mathbb{P}$ with interaction $U$, there exists a constant $\gamma > 1$ such that, for any configuration $\sigma \in \Omega$ and for any pattern $A_k \in \Omega_{C_k}$, we have*

$$(2.19) \quad \gamma^{-1} e^{-kp(U)} e^{-H(\mathscr{C}(A_k))} \leq \mathbb{P}(\mathscr{C}(A_k)) \leq \gamma e^{-kp(U)} e^{-H(\mathscr{C}(A_k))}.$$

Two other well-known properties of Gibbs measures in $d = 1$, which will be used often, are listed below.

PROPOSITION 2.10. *For the unique Gibbs measure $\mathbb{P}$ corresponding to the interaction $U \in \mathcal{U}$, there are constants $\rho < 1$ and $c > 0$, such that, for all $A_k \in \Omega_{C_k}$ and for all $\eta \in \Omega$,*

$$(2.20) \qquad \mathbb{P}(\sigma_{C_k} = A_k) \leq \rho^k$$

*and*

$$(2.21) \quad c^{-1}\mathbb{P}(\sigma_{C_k} = A_k) \leq \mathbb{P}(\sigma_{C_k} = A_k | \eta_{\mathbb{Z} \setminus C_k}) \leq \mathbb{P}(\sigma_{C_k} = A_k) c.$$

PROOF. Inequality (2.20) follows from the finite-energy property, that is, there exists $\delta > 0$ such that, for all $\sigma$,

$$0 < \delta < \mathbb{P}(\sigma_i = \alpha_i | \sigma_{\mathbb{Z} \setminus \{i\}}) < (1 - \delta).$$

This in turn follows from

$$\mathbb{P}(\sigma_i = \alpha_i | \sigma_{\mathbb{Z} \setminus \{i\}}) = \frac{\exp(-H_{\{i\}}^\sigma(\alpha_i))}{\sum_{\alpha \in \mathcal{A}} \exp(-H_{\{i\}}^\sigma(\alpha))}$$

and

$$\sup_{\sigma,\alpha_i} H^\sigma_{\{i\}}(\alpha_i) < \infty$$

by the exponential decay condition (2.8).

Therefore,

$$\mathbb{P}(\sigma_{C_k} = A_k) \leq \prod_{i \in C_k} \sup_{\sigma_{\mathbb{Z}\setminus\{i\}}} \mathbb{P}(\sigma_i = \alpha_i | \sigma_{\mathbb{Z}\setminus\{i\}}) \leq (1-\delta)^k.$$

Inequalities (2.21) are proved in [7], Proposition 8.38 and Theorem 8.39. □

2.5. *Useful lemmas.* In the proofs of our theorems we will frequently make use of the following results.

LEMMA 2.11. *For $q \geq 0$, the function $\frac{p(qU)}{q}$ is nonincreasing.*

PROOF. From the definition of $p(U)$ and $s(U)$ and from the thermodynamic relation (2.17), which is equivalent to $s = p - q\frac{dp}{dq}$, it follows immediately

$$\frac{d}{dq}\left(\frac{p(qU)}{q}\right) = -\frac{s(qU)}{q^2}.$$

The claim is then a consequence of the positivity of the entropy. □

In order to state the next lemma, we need the following notation which will be used throughout the paper.

DEFINITION 2.12. Let $a_k$ and $b_k$ be two sequences of positive numbers. Then we write

$$a_k \approx b_k,$$

if $\log(a_k) - \log(b_k)$ is a bounded sequence and

$$a_k \preceq b_k,$$

if

$$a_k \leq c_k$$

with $c_k \approx b_k$.

Note that we have that $\approx$ and $\preceq$ "behave" as ordinary equalities and inequalities and are "compatible" with usual equalities and inequalities. For example, if $a_k \preceq b_k$ and $b_k \approx c_k$, then $a_k \preceq c_k$, if $a_k \approx b_k$ and $b_k \leq c_k$, then $a_k \preceq c_k$, etc.

LEMMA 2.13. *Define*

$$\alpha = p(U) - \frac{p(2U)}{2}. \tag{2.22}$$

*We have $\alpha > 0$ and*

$$\sum_{A_k \in \Omega_{C_k}} [\mathbb{P}(\sigma_{C_k} = A_k)]^2 \approx e^{-2k\alpha}, \tag{2.23}$$

*while, for $s > 2$,*

$$\sum_{A_k \in \Omega_{C_k}} [\mathbb{P}(\sigma_{C_k} = A_k)]^s \preceq e^{-sk\alpha}. \tag{2.24}$$

PROOF. The positivity of $\alpha$ follows from Lemma 2.11. From Proposition 2.9 we obtain

$$\sum_{A_k \in \Omega_{C_k}} [\mathbb{P}(\sigma_{C_k} = A_k)]^2 \approx \sum_{A_k \in \Omega_{C_k}} e^{-2kp(U)} e^{-2H(\mathscr{C}(A_k))}$$

$$\approx e^{-2k[p(U) - p(2U)/2]} = e^{-2\alpha k}.$$

For $s > 2$, we have

$$\sum_{A_k \in \Omega_{C_k}} \mathbb{P}(\sigma_{C_k} = A_k)^s \approx \sum_{A_k \in \Omega_{C_k}} e^{-skp(U)} e^{-sH(\mathscr{C}(A_k))}$$

$$\approx e^{-sk[p(U) - p(sU)/s]} \leq e^{-s\alpha k},$$

where in the last inequality we have used the monotonicity property of Lemma 2.11. □

**3. The average number of shift matches.** We will focus on the quantity $N(\sigma, n, k)$ of Definition 2.1 and we will study how the number of shift-matchings behaves when the size of the matching, $k$, is varied as a function of the string length, $n$. It is clear that when $k = k(n)$ is very large (say, of the order of $n$), then there will be no matching of size $k$ with probability close to one, in the limit $n \to \infty$. On the other hand, if $k = k(n)$ is too small, then the number of shift-matchings will be very large with probability close to one. We want to identify a scale $k^*(n)$ such that $N(\sigma, n, k^*(n))$ will have a nontrivial distribution. Our first result concerns the average of $N(\sigma, n, k)$. Define

$$k^*(n) = \frac{\ln n}{\alpha} \tag{3.1}$$

with $\alpha$ as in (2.22). For sequences $k'(n)$ and $k(n)$, we write $k(n) \gg k'(n)$ if $k(n) - k'(n) \to \infty$ as $n \to \infty$.

Then we have the following result.

THEOREM 3.1. *Let $\{k(n)\}_{n\in\mathbb{N}}$ be a sequence of integers. Then we have the following:*

1. *If $k^*(n) \gg k(n)$, then $\lim_{n\to\infty} \mathbb{E}(N(\sigma,n,k(n))) = \infty$.*
2. *If $k(n) \gg k^*(n)$, then $\lim_{n\to\infty} \mathbb{E}(N(\sigma,n,k(n))) = 0$.*
3. *If $k(n) - k^*(n)$ is a bounded sequence, then we have*

$$(3.2) \quad 0 < \liminf_{n\to\infty} \mathbb{E}(N(\sigma,n,k(n))) \leq \limsup_{n\to\infty} \mathbb{E}(N(\sigma,n,k(n))) < \infty.$$

PROOF. We will assume (without loss of generality) that the sequence is such that

$$\lim_{n\to\infty} \frac{k(n)}{n} = 0.$$

We may rewrite $N(\sigma,n,k)$ by summing over all possible patterns of length $k$:

$$N(\sigma,n,k) = \sum_{i=0}^{n-k} \sum_{j=i+1}^{n-k} \sum_{A_k \in \Omega_{C_k}} \mathbb{1}\{(\tau_i \sigma)_{C_k} = (\tau_j \sigma)_{C_k} = A_k\}.$$

We split the above sum into two sums, one ($S_0$) corresponding to absence of overlap between $(\tau_i\sigma)_{C_k}$ and $(\tau_j\sigma)_{C_k}$ (i.e., the indices $i$ and $j$ are more than $k$ far apart) and one ($S_1$) where there is overlap:

$$S_0 = \sum_{i=0}^{n-2k} \sum_{j=i+1+k}^{n-k} \sum_{A_k \in \Omega_{C_k}} \mathbb{1}\{(\tau_i \sigma)_{C_k} = (\tau_j \sigma)_{C_k} = A_k\},$$

$$S_1 = \sum_{i=0}^{n-k} \sum_{j=i+1}^{i+k} \sum_{A_k \in \Omega_{C_k}} \mathbb{1}\{(\tau_i \sigma)_{C_k} = (\tau_j \sigma)_{C_k} = A_k\}.$$

We have of course $\mathbb{E}(N(\sigma,n,k)) = \mathbb{E}(S_0) + \mathbb{E}(S_1)$. In order to prove the first statement of the theorem, it suffices to show that $\mathbb{E}(S_0)$ diverges under the hypothesis $k^*(n) \gg k(n)$. Using translation-invariance, one has

$$\mathbb{E}(S_0) = \sum_{l=k}^{n-k}(n-k+1-l) \sum_{A_k \in \Omega_{C_k}} \mathbb{P}(\sigma_{C_k} = (\tau_l\sigma)_{C_k} = A_k)$$

$$= \sum_{l=k}^{n-k}(n-k+1-l) \sum_{A_k \in \Omega_{C_k}} \mathbb{P}(\sigma_{C_k} = A_k)\mathbb{P}((\tau_l\sigma)_{C_k} = A_k | \sigma_{C_k} = A_k).$$

Because of the mixing conditions (2.14), we have

$$(3.3) \quad \mathbb{E}(S_0) = \sum_{l=k}^{n-k}(n-k+1-l) \sum_{A_k \in \Omega_{C_k}} [\mathbb{P}(\sigma_{C_k} = A_k)]^2 + \Delta(n,k),$$

where the error $\Delta(n,k)$ is bounded by

$$|\Delta(n,k)| \leq \mathcal{O}(1) \sum_{l=k}^{n-k}(n-k+1-l) \sum_{A_k \in \Omega_{C_k}} \mathbb{P}(\sigma_{C_k} = A_k)^2 e^{-c(l-k)}.$$

Using the mixing property (2.14) and Lemma 2.13, the error can be bounded by

$$(3.4) \quad |\Delta(n,k)| \leq \mathcal{O}(1) e^{-2\alpha k} \sum_{m=0}^{n-2k} (n-2k-m+1) e^{-cm} \leq \mathcal{O}(1) e^{-2\alpha k}.$$

On the other hand, applying Lemma 2.13, we have that

$$(3.5) \quad \sum_{l=k+1}^{n-k}(n-k+1-l) \sum_{A_k} \mathbb{P}(A_k)^2 \approx (n-2k)^2 e^{-2\alpha k}.$$

Combining together (3.3), (3.4) and (3.5), we obtain

$$(3.6) \quad (n-2k)^2 e^{-2\alpha k} \preceq \mathbb{E}(N(\sigma,n,k)),$$

which proves statement 1 of the theorem.

To prove statement 2, we have to control $\mathbb{E}(S_1)$, which is the contribution to $\mathbb{E}(N(\sigma,n,k))$ due to self-overlapping cylinders. Using translation-invariance, we have

$$\mathbb{E}(S_1) = \sum_{l=1}^{k-1}(n-k+1-l) \sum_{A_k \in \Omega_{C_k}} \mathbb{P}(\sigma_{C_k} = (\tau_l \sigma)_{C_k} = A_k).$$

We further split this in two sums, namely, $\mathbb{E}(S_1) = \mathbb{E}(S_1') + \mathbb{E}(S_1'')$ with

$$(3.7) \quad \mathbb{E}(S_1') = \sum_{l=1}^{\lfloor k/2 \rfloor}(n-k+1-l) \sum_{A_k \in \Omega_{C_k}} \mathbb{P}(\sigma_{C_k} = (\tau_l \sigma)_{C_k} = A_k),$$

$$(3.8) \quad \mathbb{E}(S_1'') = \sum_{l=\lfloor k/2 \rfloor+1}^{k-1}(n-k+1-l) \sum_{A_k \in \Omega_{C_k}} \mathbb{P}(\sigma_{C_k} = (\tau_l \sigma)_{C_k} = A_k).$$

Let us consider first $\mathbb{E}(S_1'')$, that is, $\lfloor k/2 \rfloor < l < k$. In this case the overlap between $C_k$ and $\tau_l C_k$ imposes that the sum over cylinders of length $k$ can be reduced to a sum over cylinders of length $l$. In the notation of Section 2.1, we have the following inequality:

$$(3.9) \quad \begin{aligned} &\mathbb{1}(\sigma_{C_k} = (\tau_l \sigma)_{C_k} = A_k) \\ &\leq \mathbb{1}(\sigma_{C_{l+k}} = A_k(1,l) A_k(1,l) A_k(1,k-l)). \end{aligned}$$

In fact, if the pattern $A_k$ is such that the set $\{\sigma \in \Omega : \sigma_{C_k} = (\tau_l \sigma)_{C_k} = A_k\}$ is not empty, then we have equality in (3.9). Hence,

$$
\begin{aligned}
\sum_{A_k \in \Omega_k} & \mathbb{P}(\sigma_{C_k} = (\tau_l \sigma)_{C_k} = A_k) \\
&= \sum_{A_l} \sum_{B_{k-l}} \mathbb{P}(\sigma_{C_k} = A_l B_{k-l}, (\tau_l \sigma)_{C_k} = A_l B_{k-l}) \\
&\leq \sum_{A_l} \mathbb{P}(\sigma_{C_{l+k}} = A_l A_l A_l(1, k-l)) \\
&\preceq \sum_{A_l} \mathbb{P}(A_l)^2 \mathbb{P}(A_l(1, k-l)),
\end{aligned}
\tag{3.10}
$$

where in the first inequality we used the fact that contributions with $B_{k-l} \neq A_l(1, k-l)$ are zero. Therefore, using Proposition 2.10, we obtain

$$
\mathbb{E}(S_1'') \preceq \sum_{l=\lfloor k/2 \rfloor + 1}^{k} (n - k - l) \sum_{A_l} \mathbb{P}(A_l)^2 \rho^{k-l}.
$$

From this we deduce, thanks to Lemma 2.13,

$$
\begin{aligned}
\mathbb{E}(S_1'') &\preceq (n - k) \sum_{l=\lfloor k/2 \rfloor + 1}^{k} e^{-2l\alpha} \rho^{k-l} \\
&\leq (n - k) e^{-k\alpha} \sum_{l=\lfloor k/2 \rfloor + 1}^{k} \rho^{k-l} \\
&\leq (n - k) e^{-k\alpha} \sum_{x=0}^{\infty} \rho^{x} \\
&\approx (n - k) e^{-k\alpha}.
\end{aligned}
\tag{3.11}
$$

We now treat $\mathbb{E}(S_1')$, that is, the case with $1 \leq l \leq \lfloor k/2 \rfloor$. Write $k = rl + q$ with $r$ and $s$ integers, $r \geq 2$, $0 \leq q \leq l - 1$. If the set $\{\sigma : \sigma_{C_k} = (\tau_l \sigma)_{C_k} = A_k\}$ is not empty, then the pattern $A_k$ has to consist of $r + 1$ repetitions of the subpattern $A_k(1, l)$ followed by a subpattern $A_k(1, q)$, where $q$ is such that $(r+1)l + q = k + l$. Hence,

$$
\mathbb{1}(\sigma_{C_k} = (\tau_l \sigma)_{C_k} = A_k) \leq \mathbb{1}(\sigma_{C_{k+l}} = \underbrace{A_k(1, l) \cdots A_k(1, l)}_{r+1 \text{ times}} A_k(1, q)).
\tag{3.12}
$$

At this stage one could repeat the same approach as in the previous estimate for $\mathbb{E}(S_1'')$ by immediately employing Proposition 2.10. However, this approach would not work because the repeating blocks are two small. To circumvent this, we observe that in the pattern $[A_k(1, l)]^{r+1} A_k(1, q)$ there exists a piece of length $\lfloor k/2 \rfloor$ which occurs at least two times, and the remaining $l$ symbols are fixed by that piece. Therefore, using Proposition 2.10,

$$\sum_{A_k \in \Omega_k} \mathbb{P}(\sigma_{C_k} = (\tau_l \sigma)_{C_k} = A_k) \leq \sum_{B_{\lfloor k/2 \rfloor}} \mathbb{P}(B_{\lfloor k/2 \rfloor})^2 \rho^l. \tag{3.13}$$

By inserting (3.13) in (3.7) and using Lemma 2.13, we finally have

$$\mathbb{E}(S_1') \preceq (n-k)e^{-k\alpha}. \tag{3.14}$$

Combining together the estimates (3.5), (3.11) and (3.14), we obtain so far

$$\mathbb{E}(N(\sigma, n, k)) \preceq (n-k)e^{-k\alpha} + (n-2k)^2 e^{-2k\alpha} \tag{3.15}$$

from which statement 2 of the theorem follows.

Finally, combining (3.6) and (3.15) gives statement 3 of the theorem. □

**4. Second moment estimate.** In this section we will show that the random variable $N(\sigma, n, k(n))$ converges in probability to $+\infty$ in the regime where $k(n) \ll k^*(n)$, while it converges to 0 in the opposite regime $k(n) \gg k^*(n)$. Finally, if the difference $k(n) - k^*(n)$ is bounded, then we show that $N(\sigma, n, k(n))$ is tight and does not converge to zero in distribution. These results will follow as an application of the method of first moment and second moment, respectively.

THEOREM 4.1. *Let $\{k(n)\}_{n \in \mathbb{N}}$ be a sequence of integers. For every positive $m \in \mathbb{N}$:*

1. *If $k^*(n) \gg k(n)$, then $\lim_{n \to \infty} \mathbb{P}(N(\sigma, n, k(n)) \leq m) = 0$.*
2. *If $k(n) \gg k^*(n)$, then $\lim_{n \to \infty} \mathbb{P}(N(\sigma, n, k(n)) \geq m) = 0$.*
3. *If $k(n) - k^*(n)$ is bounded, then $N(\sigma, n, k(n))$ is tight and does not converge to zero in distribution. More precisely, we have that there exists a constant $C > 0$ such that*

$$\limsup_{n \to \infty} \mathbb{P}(N(\sigma, n, k(n)) > m) \leq C/m \tag{4.1}$$

*and*

$$\liminf_{n \to \infty} \mathbb{P}(N(\sigma, n, k(n)) > 0) > 0. \tag{4.2}$$

PROOF. We will assume, once more, without loss of generality that

$$\lim_{n \to \infty} \frac{k(n)}{n} = 0.$$

Statement 2 and (4.1) follow from Theorem 3.1 and the Markov inequality. To prove statement 1 and (4.2), we use the Paley–Zygmund inequality [10]

(which is an easy consequence of the Cauchy–Schwarz inequality), which gives that for all $0 \leq a \leq 1$

$$\mathbb{P}(N \geq a\mathbb{E}(N)) \geq (1-a)^2 \frac{\mathbb{E}(N)^2}{\mathbb{E}(N^2)}. \tag{4.3}$$

We fix now a sequence $k_n \uparrow \infty$ such that $k_n^* \gg k_n$. Consider the auxiliary random variable

$$\mathcal{N}_n := \sum_{i,j=0, |i-j|>2k_n}^{n-k_n} \mathbb{1}((\tau_i\sigma)_{C_{k_n}} = (\tau_j\sigma)_{C_{k_n}}). \tag{4.4}$$

Clearly, to obtain statement 1, it is sufficient that $\mathcal{N}_n$ goes to infinity with probability one. On the other hand, using the first moment computations of the previous section, we have

$$\mathbb{E}(\mathcal{N}_n) \approx n^2 e^{-2\alpha k_n}. \tag{4.5}$$

So, in order to use the Paley–Zygmund inequality, it is sufficient to show that

$$\mathbb{E}(\mathcal{N}_n^2) \preceq \xi_n^4, \tag{4.6}$$

where we introduced the notation

$$\xi_n := n e^{-\alpha k_n}. \tag{4.7}$$

Remark that $\xi_n \to \infty$ for our choice of $k_n$ (as in statement 1).

Indeed, if we have (4.6) in the regime $k^*(n) \gg k(n)$, then the ratio

$$\frac{\mathbb{E}(\mathcal{N}^2)}{(\mathbb{E}(\mathcal{N}))^2}$$

remains bounded from above as $n \to \infty$, and hence, using (4.3), $\mathcal{N}_n$ diverges with probability at least $\delta > 0$. Therefore, in that case, by ergodicity, $N(\sigma, n, k_n) \geq \mathcal{N}_n$ goes to infinity with probability one, since the set of $\sigma$'s such that $N(\sigma, n, k_n)$ goes to infinity is translation-invariant, and hence has measure zero or one.

To see how statement (4.2) follows from (4.6) in the regime where $k(n) - k^*(n)$ is bounded, use the (more classical) second moment inequality

$$\mathbb{P}(\mathcal{N} > 0) \geq \frac{(\mathbb{E}(\mathcal{N}))^2}{\mathbb{E}(\mathcal{N}^2)}$$

combined with

$$N(\sigma, n, k(n)) \geq \mathcal{N}.$$

We now proceed with the proof of (4.6). We have

$$\mathbb{E}(\mathcal{N}_n^2) = \sum_{i,j,r,s,|i-j|>2k_n,|r-s|>2k_n} \sum_{A_{k_n},B_{k_n}} \mathbb{P}((A_{k_n})_i(A_{k_n})_j(B_{k_n})_r(B_{k_n})_s),$$
(4.8)

where we use the abbreviate notation $(A_{k_n})_i$ for the event $(\tau_i\sigma)_{C_{k_n}} = A_{k_n}$. Similarly, if we have a word of length $l$, say, consisting of p symbols of $A_p$ followed by $l-p$ symbols of $B_{l-p}$, we write $(A_p B_{l-p})_i$ for the event that this word appears at location $i$, that is, the event $(\tau_i\sigma)_{C_l} = A_p B_{l-p}$.

The sum in the right-hand side of (4.8) will be split into different sums, according to the amount of overlap in the set of indices $\{i,j,r,s\}$. By this we mean the following: we say that there is *overlap* between two indices $i,j$ if $|i-j| < k_n$. The number of overlaps of a set of indices $\{i,j,r,s\}$ is denoted by $\theta(i,j,r,s)$ and is the number of unordered pairs of indices which have overlap. Since we restrict in the sum (4.8) to $|i-j| > 2k_n, |r-s| > 2k_n$, it follows from the triangular inequality that in that case $\theta(i,j,r,s) \leq 2$. Therefore, we split the sum into three cases

(4.9) $$\sum_{i,j,r,s,|i-j|>2k_n,|r-s|>2k_n} \sum_{A_{k_n},B_{k_n}} \mathbb{P}((A_{k_n})_i(A_{k_n})_j(B_{k_n})_r(B_{k_n})_s)$$
$$= S_0 + S_1 + S_2,$$

where

(4.10) $$S_p = \sum_{(i,j,r,s)\in K_{k,p}} \sum_{A,B} \mathbb{P}((A_{k_n})_i(A_{k_n})_j(B_{k_n})_r(B_{k_n})_s),$$

where we abbreviated

(4.11) $$K_{k_n,p} = \{(i,j,r,s) : |i-j| > 2k_n, |r-s| > 2k_n, \theta(i,j,r,s) = p\}$$

to be the set of indices such that the overlap is $p$.

1. *Zero overlap*: $S_0$.

We use Lemma 2.13, and notation (4.7):

(4.12) $$S_0 \preceq \sum_{i,j,r,s} \sum_{A_{k_n},B_{k_n}} \mathbb{P}(A_{k_n})^2 \mathbb{P}(B_{k_n})^2 \preceq \xi_n^4.$$

2. *One overlap*: $S_1$.

We treat the case $|i-r| < k_n$, $i < r < j < s$. The other cases are treated in exactly the same way. Put $A_{k_n} = [a_1, a_2, \ldots, a_{k_n}]$, $B_{k_n} = [b_1, b_2, \ldots, b_{k_n}]$. The intersection $(A_{k_n})_i \cap (B_{k_n})_r$ is nonempty if and only if $a_r = b_1$, $a_{r+1} = b_2, \ldots, a_{k_n} = b_{k_n-r+1}$, that is, the last $k_n - r + 1$ symbols of $A_{k_n}$ are equal to the first $k_n - r + 1$ symbols of $B_{k_n}$.

Therefore, we obtain that the sum over the patterns $A_{k_n}, B_{k_n}$ in $S_1$ equals

$$\sum_{A_{k_n},B_{k_n}} \mathbb{P}((A_{k_n})_i(A_{k_n})_j(B_{k_n})_r(B_{k_n})_s)$$

$$= \sum_{A_{k_n},B_{k_n}} \mathbb{P}((A_{k_n}B_{k_n}(k_n-r,k_n))_i(A_{k_n})_j(A_{k_n}(r,k_n)B_{k_n}(k_n-r,k_n))_s)$$

(4.13)

$$\preceq \sum_{A_{k_n},B_{k_n}} \mathbb{P}(A_{k_n}(r,k_n))^3 \mathbb{P}(A_{k_n}(1,r-1))^2 \mathbb{P}(B_{k_n}(k_n-r,k_n))^2$$

$$\preceq e^{-3(k_n-r)\alpha} e^{-2r\alpha} e^{-2r\alpha}.$$

Summing over the indices $(i,j,r,s) \in K(k_n,1)$ then gives

(4.14) $$S_1 \preceq n^3 e^{-3\alpha k_n} \sum_{r \leq k_n} e^{-r\alpha} \preceq \xi_n^3.$$

3. *Two overlaps*: $S_2$.

We treat the case $i < r < j < s$ and $r - i < k_n, s - j < k_n$. Other cases are treated in the same way. Put $l_1 := i + k_n - r + 1, p_1 = j + k_n - s + 1$. We suppose $l_1 > p_1$. Then the last $l_1$ symbols of $A_{k_n}$ have to equal the first $l_1$ symbols of $B_{k_n}$, otherwise the intersection $(A_{k_n})_i(A_{k_n})_j(B_{k_n})_r(B_{k_n})_s$ is empty. Therefore, we obtain that the sum over the patterns $A_{k_n}, B_{k_n}$ in $S_2$ equals

$$\sum_{A_{k_n},B_{k_n}} \mathbb{P}((A_{k_n})_i(A_{k_n})_j(B_{k_n})_r(B_{k_n})_s)$$

$$= \sum_{A_{k_n},B_{k_n}} \mathbb{P}((A_{k_n}B_{k_n-l_1})_i(A_{k_n}B_{k_n-p_1})_j)$$

(4.15)

$$\preceq \sum_{A_{k_n},B_{k_n}} \mathbb{P}(A_{k_n})^2 \mathbb{P}(B_{k_n-l_1})^2 \rho^{l_1-p_1}$$

$$\preceq e^{-2k\alpha} e^{-2(k-l_1)\alpha} \rho^{l_1-p_1}.$$

Summing over the indices in $K(k,2)$ then gives

(4.16) $$S_2 \preceq n^2 e^{-2k_n\alpha} \sum_{l_1 < k_n} e^{-2\alpha(k_n-l_1)} \sum_{p_1 < l_1} \rho^{l_1-p_1} \preceq \xi_n^2.$$

Using the bounds (4.12), (4.14) and (4.16) in (4.8) and (4.9), we deduce (4.6) and then, as explained below, statement 1 of the theorem follows from the Paley–Zygmund inequality. This completes the proof. □

The following result relates Theorem 4.1 and the behavior of the maximal shift-matching, and is the analogue of Theorem 1 in [6] (which is, however, convergence almost surely for more general comparison of sequences based on scores, but for independent sequences).

PROPOSITION 4.2. *Let $M(\sigma,n)$ be defined as in Definition 2.2. Recall*

$$\alpha = p(U) - \frac{p(2U)}{2}.$$

*Then we have that*

$$\frac{M(\sigma,n)}{n} \to \alpha,$$

*where the convergence is in probability.*

PROOF. Use the relations of Proposition 2.4. We have

$$\mathbb{P}\left(\frac{M(\sigma,n)}{\alpha \log n} \geq (1+\varepsilon)\right) \leq \mathbb{P}(N(\sigma,n,\lfloor \alpha(1+\varepsilon)\log n \rfloor) \geq 1)$$

and

$$\mathbb{P}\left(\frac{M(\sigma,n)}{\alpha \log n} < (1-\varepsilon)\right) \leq \mathbb{P}(N(\sigma,n,\lceil \alpha(1-\varepsilon)\log n \rceil) = 0).$$

So the result follows from Theorem 4.1. □

**5. Two independent strings.** In this section we study the number of matches with shift when two *independent* sequences $\sigma$ and $\eta$ are considered. The marginal distributions of $\sigma$ and $\eta$ are denoted with $\mathbb{P}$ and $\mathbb{Q}$, which are chosen to be Gibbs measure with exponentially decaying translation-invariant interactions $U(X,\sigma)$ and $V(X,\eta)$, respectively. We assume the two strings belong to the same alphabet $\mathcal{A}$. In analogy with the case of one string, we give the following definition.

DEFINITION 5.1 (Number of shift-matches for 2 strings). For every couple of configurations $\sigma,\eta \in \Omega \times \Omega$ and for every $n \in \mathbb{N}$, $k \in \mathbb{N}$, with $k < n$, we define the number of matches with shift of length $k$ as

$$(5.1) \qquad N(\sigma,\eta,n,k) = \sum_{i=0}^{n-k} \sum_{j=0, j\neq i}^{n-k} \mathbb{1}\{(\tau_i\sigma)_{C_k} = (\tau_j\eta)_{C_k}\}.$$

Of course, in the case $\sigma = \eta$ we recover (up to a factor 2) the previous Definition 2.1, that is, $N(\sigma,\sigma,n,k) = 2N(\sigma,n,k)$.

5.1. *Identical marginal distribution.* We treat here the case $\mathbb{Q} = \mathbb{P}$, that is, the two sequences $\sigma$ and $\eta$ are chosen independently from the same Gibbs distribution $\mathbb{P}$ with interaction $U(X,\sigma)$. Then the results of the previous section are generalized as follows.

THEOREM 5.2. *Let $\{k(n)\}_{n\in\mathbb{N}}$ be a sequence of integers:*

1. If $k^*(n) \gg k(n)$, then $\lim_{n\to\infty} \mathbb{E}_{\mathbb{P}\otimes\mathbb{P}}[N(\sigma,\eta,n,k(n))] = \infty$.
2. If $k^*(n) \ll k(n)$, then $\lim_{n\to\infty} \mathbb{E}_{\mathbb{P}\otimes\mathbb{P}}[N(\sigma,\eta,n,k(n))] = 0$.
3. If $k(n) - k^*(n)$ is a bounded sequence, then we have

$$
\begin{aligned}
0 &< \liminf_{n\to\infty} \mathbb{E}_{\mathbb{P}\otimes\mathbb{P}}(N(\sigma,\eta,n,k(n))) \\
&\leq \limsup_{n\to\infty} \mathbb{E}_{\mathbb{P}\otimes\mathbb{P}}(N(\sigma,\eta,n,k(n))) < \infty.
\end{aligned}
\tag{5.2}
$$

PROOF. Because of independence, we immediately have

$$
\begin{aligned}
\mathbb{E}_{\mathbb{P}\otimes\mathbb{P}}&[N(\sigma,\eta,n,k)] \\
&= \sum_{i\neq j=0}^{n-k} \sum_{A_k \in \Omega_k} \mathbb{P}((\tau_i\sigma)_{C_k} = A_k)\mathbb{P}((\tau_j\eta)_{C_k} = A_k) \\
&= (n-k)^2 \sum_{A_k \in \Omega_{C_k}} \mathbb{P}(A_k)^2 \\
&\approx (n-k)^2 e^{-2k\alpha}.
\end{aligned}
\tag{5.3}
$$
$\square$

THEOREM 5.3. *Let $\{k(n)\}_{n\in\mathbf{N}}$ be a sequence of integers. For every positive $m \in \mathbb{N}$:*

1. *If $k^*(n) \gg k(n)$, then $\lim_{n\to\infty} \mathbb{P}\otimes\mathbb{P}[N(\sigma,\eta,n,k(n)) \leq \varepsilon] = 0$.*
2. *If $k^*(n) \ll k(n)$, then $\lim_{n\to\infty} \mathbb{P}\otimes\mathbb{P}[N(\sigma,\eta,n,k(n)) \geq \varepsilon] = 0$.*
3. *If $k(n) - k^*(n)$ is bounded, then $N(\sigma,\eta,n,k(n))$ is tight and does not converge to zero in distribution. More precisely, we have that there exists a constant $C > 0$ such that*

$$\limsup_{n\to\infty} \mathbb{P}\otimes\mathbb{P}(N(\sigma,\eta,n,k(n)) > m) \leq C/m \tag{5.4}$$

*and*

$$\liminf_{n\to\infty} \mathbb{P}\otimes\mathbb{P}(N(\sigma,\eta,n,k(n)) > 0) > 0. \tag{5.5}$$

PROOF. The strategy of the proof is as in Theorem 4.1. Thus, we need to control the second moment to show that $\mathbb{E}(N^2) \approx (\mathbb{E}(N))^2$. We start from

$$
\begin{aligned}
\mathbb{E}_{\mathbb{P}\otimes\mathbb{P}}&(N^2(\sigma,\eta,n,k)) \\
&= \sum_{i_1,j_1,i_2,j_2=1}^{n-k} \sum_{A_k,B_k \in \Omega_k} \mathbb{P}((\tau_{i_1}\sigma)_{C_k} = A_k, (\tau_{i_2}\sigma)_{C_k} = B_k) \\
&\qquad \times \mathbb{P}((\tau_{j_1}\eta)_{C_k} = A_k, (\tau_{j_2}\eta)_{C_k} = B_k).
\end{aligned}
\tag{5.6}
$$

Using translation-invariance and defining new indices $l_1 = i_2 - i_1$ and $l_2 = j_2 - j_1$, we have

$$\mathbb{E}_{\mathbb{P}\otimes\mathbb{P}}(N^2(\sigma,\eta,n,k))$$
$$= \sum_{A_k,B_k\in\Omega_k} \left( \sum_{l_1=1}^{n-k} (n-k+1-l_1)\mathbb{P}(\sigma_{C_k} = A_k, (\tau_{l_1}\sigma)_{C_k} = B_k) \right.$$
$$\left. \times \sum_{l_2=1}^{n-k} (n-k+1-l_2)\mathbb{P}(\eta_{C_k} = A_k, (\tau_{l_2}\eta)_{C_k} = B_k) \right).$$

We have to distinguish three kinds of contributions in the previous sums:

1. Zero overlap, that is, $l_1 > k, l_2 > k$. Then

$$\sum_{A_k,B_k\in\Omega_k} \left( \sum_{l_1=k+1}^{n-k} (n-k+1-l_1)\mathbb{P}(\sigma_{C_k} = A_k, (\tau_{l_1}\sigma)_{C_k} = B_k) \right.$$
$$\left. \times \sum_{l_2=k+1}^{n-k} (n-k+1-l_2)\mathbb{P}(\eta_{C_k} = A_k, (\tau_{l_2}\eta)_{C_k} = B_k) \right)$$

(5.7)
$$\approx (n-k)^4 \sum_{A_k,B_k\in\Omega_k} \mathbb{P}(A_k)^2 \mathbb{P}(B_k)^2$$
$$\approx (n-k)^4 e^{-4k\alpha}.$$

2. One overlap. We treat the case $l_1 \leq k$ and $l_2 > k$ (other cases are treated similarly). We have

$$\sum_{A_k,B_k\in\Omega_k} \left( \sum_{l_1=1}^{k} (n-k+1-l_1)\mathbb{P}(\sigma_{C_k} = A_k, (\tau_{l_1}\sigma)_{C_k} = B_k) \right.$$
$$\left. \times \sum_{l_2=k+1}^{n-k} (n-k+1-l_2)\mathbb{P}(\eta_{C_k} = A_k, (\tau_{l_2}\eta)_{C_k} = B_k) \right)$$

$$\approx (n-k)^3 \sum_{l_1=1}^{k} \sum_{D_{l_1},E_{k-l_1},F_{l_1}} \mathbb{P}(D_{l_1}E_{k-l_1}F_{l_1})\mathbb{P}(D_{l_1}E_{k-l_1})\mathbb{P}(E_{k-l_1}F_{l_1})$$

(5.8)
$$\approx (n-k)^3 \sum_{l_1=1}^{k} \sum_{D_{l_1},E_{k-l_1},F_{l_1}} \mathbb{P}(D_{l_1})^2 \mathbb{P}(E_{k-l_1})^3 \mathbb{P}(F_{l_1})^2$$

$$\preceq (n-k)^3 \sum_{l_1=1}^{k} e^{-2l_1\alpha} e^{-2l_1\alpha} e^{-3(k-l_1)\alpha}$$

$$\leq (n-k)^3 e^{-3k\alpha}.$$

3. Two overlaps. We treat the case $l_1 < l_2 \leq k$ (other cases are treated similarly). We have

$$\sum_{A_k, B_k \in \Omega_k} \left( \sum_{l_1=1}^{k} (n-k+1-l_1) \mathbb{P}(\sigma_{C_k} = A_k, (\tau_{l_1}\sigma)_{C_k} = B_k) \right.$$
$$\left. \times \sum_{l_2=1}^{k} (n-k+1-l_2) \mathbb{P}(\eta_{C_k} = A_k, (\tau_{l_2}\eta)_{C_k} = B_k) \right)$$

$$\approx (n-k)^2 \sum_{l_1, l_2=1}^{k} \sum_{D_{l_1}, E_{l_2-l_1}, F_{k-l_2}, G_{l_1}, H_{l_2-l_1}} \mathbb{P}(D_{l_1} E_{l_2-l_1} F_{k-l_2} G_{l_1})$$
$$\times \mathbb{P}(D_{l_1} E_{l_2-l_1} F_{k-l_2} G_{l_1} H_{l_2-l_1})$$

(5.9)

$$\approx (n-k)^2 \sum_{l_1, l_2=1}^{k} \sum_{D_{l_1}} \mathbb{P}(D_{l_1})^2 \sum_{E_{l_2-l_1}} \mathbb{P}(E_{l_2-l_1})^2 \sum_{F_{k-l_2}} \mathbb{P}(F_{k-l_2})^2$$
$$\times \sum_{G_{l_1}} \mathbb{P}(G_{l_1})^2 \sum_{H_{l_2-l_1}} \mathbb{P}(H_{l_2-l_1})$$

$$\preceq (n-k)^2 \sum_{l_1, l_2=1}^{k} e^{-2l_1\alpha} e^{-2(l_2-l_1)\alpha} e^{-2(k-l_2)\alpha} e^{-2l_1\alpha}$$

$$\leq (n-k)^2 e^{-2k\alpha}.$$

Combining together (5.7), (5.8) and (5.9) and similar expression for other cases with one and two overlaps, we obtain the second moment condition $\mathbb{E}(N^2) \preceq (\mathbb{E}(N))^2$. $\square$

5.2. *Different marginal distributions.* In the case $\mathbb{P} \neq \mathbb{Q}$, the first moment is controlled in an analogous way, but the second moment analysis is different, and, in fact, as we will show in an example, it can happen for some scale $k_n \to \infty$ that:

1. $\mathbb{E}_{\mathbb{P} \otimes \mathbb{Q}}(N(\sigma, \eta, n, k_n)) \to \infty$ as $n \to \infty$,
2. $\mathbb{P} \otimes \mathbb{Q}(N(\sigma, \eta, n, k_n) = 0) > e^{-\delta}$ for some $\delta > 0$ independent of $n$.

This means that in order to decide whether $N(\sigma, \eta, n, k_n)$ goes to infinity $\mathbb{P} \otimes \mathbb{Q}$ almost surely, it is not sufficient to have $\mathbb{E}_{\mathbb{P} \otimes \mathbb{Q}}(N(\sigma, \eta, n, k_n)) \to \infty$.

We start with the case $\mathbb{P}$ and $\mathbb{Q}$ Gibbs measures with potentials $U, V$, respectively, and define

(5.10) $$\tilde{\alpha} = \tfrac{1}{2}p(U) + \tfrac{1}{2}p(V) - \tfrac{1}{2}p(U+V) > 0$$

and

(5.11) $$\tilde{k}^* = \frac{\log n}{\tilde{\alpha}},$$

then we have the following:

THEOREM 5.4. *Let $\{k(n)\}_{n\in\mathbb{N}}$ be a sequence of integers.*

1. *If $\tilde{k}^*(n) \gg k(n)$, then $\lim_{n\to\infty} \mathbb{E}_{\mathbb{P}\otimes\mathbb{Q}}(N(\sigma,\eta,n,k(n))) = \infty$.*
2. *If $\tilde{k}^*(n) \ll k(n)$, then $\lim_{n\to\infty} \mathbb{E}_{\mathbb{P}\otimes\mathbb{Q}}(N(\sigma,\eta,n,k(n))) = 0$.*
3. *If $k(n) - \tilde{k}^*(n)$ is a bounded sequence, then we have*

$$
\begin{aligned}
0 &< \liminf_{n\to\infty} \mathbb{E}_{\mathbb{P}\otimes\mathbb{Q}}(N(\sigma,\eta,n,k(n))) \\
&\leq \limsup_{n\to\infty} \mathbb{E}_{\mathbb{P}\otimes\mathbb{Q}}(N(\sigma,\eta,n,k(n))) < \infty.
\end{aligned}
\tag{5.12}
$$

PROOF. Start by rewriting

$$N(\sigma,\eta,n,k) = \sum_{i=0}^{n-k} \sum_{j=0, j\neq i}^{n-k} \sum_{A_k \in \Omega_k} \mathbb{1}\{(\tau_i\sigma)_{C_k} = A_k, (\tau_j\eta)_{C_k} = A_k\}.$$

Taking into account the independence of the measures $\mathbb{P}$ and $\mathbb{Q}$, we obtain

$$
\begin{aligned}
&\mathbb{E}_{\mathbb{P}\otimes\mathbb{Q}}(N(\sigma,\eta,n,k)) \\
&= \sum_{i\neq j=0}^{n-k} \sum_{A_k \in \Omega_k} \mathbb{P}((\tau_i\sigma)_{C_k} = A_k)\mathbb{Q}((\tau_j\eta)_{C_k} = A_k) \\
&\approx (n-k)^2 \sum_{A_k \in \Omega_{C_k}} e^{-kp(U)} e^{-kH_U(\mathscr{C}(A_k))} e^{-kp(V)} e^{-kH_V(\mathscr{C}(A_k))} \\
&\approx (n-k)^2 e^{-k[p(U)+p(V)-p(U+V)]} \\
&= (n-k)^2 e^{-2k\tilde{\alpha}},
\end{aligned}
\tag{5.13}
$$

where in the second line we made use of translation-invariance and Proposition 2.9. □

In case 1 of Theorem 5.4, we will not in general be able to conclude that $N(\sigma,\eta,n,k(n))$ goes to infinity almost surely as $n \to \infty$. Indeed, if we compute the second moment, we find terms analogous to the case $\mathbb{P} = \mathbb{Q}$, of which now we have to take the $\mathbb{P} \otimes \mathbb{Q}$ expectation. In particular, the one overlap contribution will contain a term of the order

$$(n-k)^3 \sum_{E_k} \mathbb{P}(E_k)\mathbb{Q}(E_k)^2.$$

If $\mathbb{P} \neq \mathbb{Q}$, this term may however not be dominated by $n^4 e^{-4k\tilde{\alpha}}$. Indeed, the inequality

$$\sum_{E_k} \mathbb{P}(E_k)\mathbb{Q}(E_k)^2 \leq \left(\sum_{E_k} \mathbb{P}(E_k)\mathbb{Q}(E_k)\right)^{3/2}$$

is not valid in general. In particular, if $\mathbb{P}$ gives uniform measure to cylinders $E_k$ and $\mathbb{Q}$ concentrates on one particular cylinder, then this inequality will be violated.

As an example, inspired by this, we choose $\mathbb{P}$ to be a Gibbs measure with potential $U$, and $\mathbb{Q} = \delta_a$, where $\delta_a$ denotes the Dirac measure concentrating on the configuration $\eta(x) = a$ for all $x \in \mathbb{Z}$ (which is strictly speaking not a Gibbs measure, but a limit of Gibbs measures). In that case $\mathbb{P} \otimes \mathbb{Q}$ almost surely,

$$N(\sigma, \eta, n, k(n)) = n \sum_{i=1}^{n-k} \mathbb{1}((\tau_i \sigma)_{C_k} = [a]_k),$$

where $[a]_k$ denotes a block of $k$ successive $a$'s. Therefore,

$$\mathbb{P} \otimes \mathbb{Q}(N(\sigma, \eta, n, k(n)) = 0) = \mathbb{P}(\Theta_{[a]_k}(\sigma) \geq n - k),$$

where

$$\Theta_{[a]_k}(\sigma) = \inf\{j > 0 : \sigma_j = a, \sigma_{j+1} = a, \ldots, \sigma_{j+k-1} = a\}$$

is the hitting time of the pattern $[a]_k$ in the configuration $\sigma$. For this hitting time we have the exponential law [1, 2] which gives

$$\mathbb{P}(\Theta_{[a]_k}(\sigma) \geq n) \geq e^{-\lambda \mathbb{P}([a]_k) n}$$

with $\lambda$ a positive constant not depending on $n$. Now we choose the scale $k_n$ such that the first moment of $N(\sigma, \eta, n, k(n))$ diverges as $n \to \infty$, that is, such that

$$n^2 \mathbb{P}([a]_{k_n}) \to \infty.$$

Furthermore, we impose that

$$\mathbb{P}([a]_{k_n}) n \leq \delta$$

for all $n$. In that case

$$\mathbb{P}(\Theta_{[a]_{k_n}}(\sigma) \geq n) \geq e^{-\lambda \mathbb{P}([a]_{k_n})) n} \geq e^{-\lambda \delta},$$

which implies $N(\sigma, \eta, n, k_n)$ does not go to infinity $\mathbb{P} \otimes \mathbb{Q}$ almost surely.

**Acknowledgment.** We thank the anonymous referee for helpful remarks and a careful reading.

# REFERENCES

[1] ABADI, M. (2001). Exponential approximation for hitting times in mixing processes. *Math. Phys. Electron. J.* **7** 19. MR1871384
[2] ABADI, M., CHAZOTTES, J.-R., REDIG, F. and VERBITSKIY, E. (2004). Exponential distribution for the occurrence of rare patterns in Gibbsian random fields. *Comm. Math. Phys.* **246** 269–294. MR2048558
[3] BOWEN, R. (2008). *Equilibrium States and the Ergodic Theory of Anosov Diffeomorphisms*, revised ed. *Lecture Notes in Mathematics* **470**. Springer, Berlin. MR2423393
[4] COLLET, P., GALVES, A. and SCHMITT, B. (1999). Repetition times for Gibbsian sources. *Nonlinearity* **12** 1225–1237. MR1709841
[5] DEMBO, A., KARLIN, S. and ZEITOUNI, O. (1994). Limit distribution of maximal non-aligned two-sequence segmental score. *Ann. Probab.* **22** 2022–2039. MR1331214
[6] DEMBO, A., KARLIN, S. and ZEITOUNI, O. (1994). Critical phenomena for sequence matching with scoring. *Ann. Probab.* **22** 1993–2021. MR1331213
[7] GEORGII, H.-O. (1988). *Gibbs Measures and Phase Transitions. de Gruyter Studies in Mathematics* **9**. de Gruyter, Berlin. MR956646
[8] GUYON, X. (1995). *Random Fields on a Network: Modeling, Statistics, and Applications*. Springer, New York. MR1344683
[9] HANSEN, N. R. (2006). Local alignment of Markov chains. *Ann. Appl. Probab.* **16** 1262–1296. MR2260063
[10] PALEY, R. and ZYGMUND, A. (1932). A note on analytic functions in the unit circle. *Proc. Camb. Phil. Soc.* **28** 266–272.
[11] RUELLE, D. (1978). *Thermodynamic Formalism: The Mathematical Structures of Classical Equilibrium Statistical Mechanics. Encyclopedia of Mathematics and Its Applications* **5**. Addison-Wesley, Reading, MA. MR511655

P. COLLET
CENTRE DE PHYSIQUE THÉORIQUE
CNRS UMR 7644
91128 PALAISEAU CEDEX
FRANCE
E-MAIL: collet@cpht.polytechnique.fr

C. GIARDINA
DEPARTMENT OF MATHEMATICS
 AND COMPUTER SCIENCE
EINDHOVEN UNIVERSITY
P.O. BOX 513—5600 MB EINDHOVEN
THE NETHERLANDS
E-MAIL: c.giardina@tue.nl

F. REDIG
MATHEMATISCH INSTITUUT
UNIVERSITEIT LEIDEN
NIELS BOHRWEG 1
2333 CA LEIDEN
THE NETHERLANDS
E-MAIL: redig@math.leidenuniv.nl